\documentclass[11pt]{amsart}
\usepackage{amsmath,amsfonts,amsthm,amssymb,amscd}
\def\classification#1{\def\@class{#1}}
\classification{\null}

\DeclareFontFamily{OT1}{rsfs}{}
\DeclareFontShape{OT1}{rsfs}{n}{it}{<-> rsfs10}{}
\DeclareMathAlphabet{\mathscr}{OT1}{rsfs}{n}{it}

\newcommand{\R}{{\mathbb R}}

\newcommand{\C}{\mathbb{C}}

\newtheorem{theorem}{Theorem}

\theoremstyle{remark}
\newtheorem{remark}[theorem]{Remark}


\title{On the number of classes of triangles determined by $N$ points in $\R^2$}

\author{Misha Rudnev}
\address{Misha Rudnev: Department of Mathematics, University of Bristol, Bristol BS8 1TW, UK}
\email{m.rudnev@bristol.ac.uk}

\begin{document}

\begin{abstract} Let $P$ be a set  of $N$ points in the Euclidean plane, where a positive proportion of points lies off a single straight line. This note points out two facts concerning the number of equivalence classes of triangles that $P$ determines, namely that (i) $P$ determines $\Omega(N^2)$ different equivalence classes of congruent triangles, and (ii) $P$ determines $\Omega(\frac{N^2}{\log N})$ different equivalence classes of similar triangles. The first fact follows from the recent theorem by Guth-Katz on point-line incidences in $\R^3$ (\cite{GK}). The second one, perhaps not so well known, is due to Solymosi and Tardos (\cite{SoTa}).

\end{abstract}

\maketitle

\section{Introduction and statement of results}

\vskip.125in
After the original version of this note was released, the author became aware of the  2007 work of Solymosi and Tardos (\cite{SoTa}) which, modulo an application of the Cauchy-Schwarz inequality subsumes the claim (ii) of the main Theorem \ref{mainth} herein. It appears nonetheless reasonable to retain the current revised version available through the {\sf arXiv} in order to juxtapose how the two closely related geometric problems about triangles in the Euclidean plane get analysed  in terms of two different group actions, as well as to bring attention to the results and questions raised in (\cite{ST}), which apart from the action of the isometries and linear complex transformations on $\R^2$ deals with the action of the M\"obius group on the Riemann sphere.

\medskip
The basic question of extremal geometric combinatorics is to give some universal statistics for the number of distinct types of certain geometric configurations that a sufficiently large point set can determine. See, for example \cite{BMP} and the references contained therein.
Today there is a  new paradigm, often credited to the work of Elekes and Sharir (\cite{ES}), which has recently resulted in a considerable progress in the analysis of such problems. It contributed an essential building block to the recent achievement by Guth and Katz (\cite{GK}) who settled the long-standing Erd\H os distance conjecture (\cite{E}) by proving that a planar set $E$  of $N$ points determines $\Omega(\frac{N}{\log N})$ distinct distances, i.e. distinct congruence classes of line segments. Shortly thereafter results of similar flavour were obtained in \cite{IRR}, and \cite{RR}, proving similar estimates for the number of dot products and Minkowski distances, respectively, and deriving form them new sum-product type estimates. However, the idea in question had had applications prior to \cite{ES} as well, in particular in \cite{SoTa}.

The idea is to lift the problem into higher dimensions by considering symmetries acting on pairs of objects within each equivalence class and try to do the statistics on the number of symmetries. In the recent applications of this idea which originated in \cite{GK}, the original problem in the plane would be adequately represented by an incidence problem between a set of straight lines and points in $\R^3$. To this end Guth and Katz (\cite{GK}) proved the following incidence theorem. (The next two theorems are given a formulation adapted for the purposes of this note).

\begin{theorem}\label{theorem:GKmain} Let $L$ be a set of $N^2$ lines in  $\mathbb{R}^3$, such that
 no more than $N$ lines are co-planar or concurrent. Then, for $k\geq 3$, the number of points where $k$ or more lines intersect is $O(\frac{N^3}{k^2}).$
\end{theorem}
Not only have  Guth and Katz proven Theorem \ref{theorem:GKmain} (its special case $k=3$ was proven prior to that by Elekes, Kaplan, and Sharir, \cite{EKS}), but they have thereby set forth a new method for proving point-line incidence theorems in $\R^d$, $d\geq 2$, based on the space decomposition provided by the polynomial version of the Ham Sandwich theorem of Stone and Tukey (\cite{StT}) attributed to Gromov (\cite{G}). The realisation of this fact is due (among others) to
Kaplan, Matou$\check{\rm s}$ek, and Sharir (\cite{KMS}),
 Solymosi and Tao (\cite{SoT}) and Zahl (\cite{Z}, \cite{Z1}), who have demonstrated that this method potentially enables one to prove the whole family of incidence theorems, not only in $\R^d$, but $\C^d$ as well.

In particular, this concerns the classical Szemer\'edi-Trotter theorem.
\begin{theorem}\label{theorem:ST} Let $L$ be a set of $N^2$ lines in  the plane, such that
 no more than $N$ lines are concurrent. Then, for $k\geq 2$, the number of points where $k$ or more lines intersect is $O(\frac{N^4}{k^3}).$
\end{theorem}
Note that the bound of Theorem \ref{theorem:GKmain} is better than that of Theorem \ref{theorem:ST}, hence the non-planarity constraint of the former theorem.\footnote{The constraint $k\geq 3$ rather than $2$ in the former theorem is to avoid another non-degeneracy condition. Naturally, the number of all kinds of geometric obstacles increases dramatically with $d$, see e.g. the recent work of Zahl (\cite{Z1}) in four dimensions.}

For this note's modest purpose, one shall only need the above two theorems, but as for Theorem \ref{theorem:ST}, it will be used over $\C^2$, rather than $\R^2$. The proof that the Szemer\'edi-Trotter theorem holds in the complex plane is originally due to T\'oth (\cite{T}). A more modern, in light of the current discussion, proof has been recently given by Zahl (\cite{Z1}). As a matter of fact, for the specific purpose of (ii) one can use a particular, and easier version of the complex Szemer\'edi-Trotter theorem, see \cite{SoTa}.

\medskip
The main result pointed out in this note is as follows.
\begin{theorem}\label{mainth} Let $P$ be a set  of $N\gg 1$ points in $\R^2$, with no single straight line containing more than $\frac{N}{2}$ points. Then

(i) $P$ determines $\Omega(N^2)$ distinct equivalence types of congruent triangles,

(ii) $P$ determines $\Omega(\frac{N^2}{\log N})$ distinct equivalence types of similar, i.e. homothetic, triangles.
\end{theorem}

\begin{remark} The example when $P$ is a square subset of the integer lattice shows that the bound in (i) is sharp. The author is prone to believe, although the number of equivalence classes by similarity is generally smaller than by congruence, that the logarithmic term in the bound (ii) is the method's artefact, however has only indirect evidence to it.\end{remark}

The statement (i) of the theorem is a straightforward corollary of the work of Guth and Katz (\cite{GK}). The statement (ii) follows, after an application of the Cauchy-Schwarz inequality from Theorem 3 of Solymosi and Tardos (\cite{ST}). Since there is basically a one-to-one set-up for the proof of both parts, it is shown in some detail in the remaining section. 

\subsubsection*{Notation} Above and throughout, $|\cdot|$ denotes cardinality of a finite set. The standard notation $X\ll Y,$ or equivalently $X=O(Y)$, means that there exists a universal constant $C>0$, such that $X\leq CY$. Conversely, $X\gg Y,$ or equivalently $X=\Omega(Y)$, means that there exists a universal constant $c>0$, such that $X\geq cY$. We write $X\approx Y$ if $X\ll Y$ and $Y\ll X$. The number $N\gg 1$ is viewed as an asymptotic parameter.

\section{Proof of Theorem \ref{mainth}}

\begin{proof} Let $T(P)$ be the set of all triangles determined by triples of distinct points of $P$. The non-collinearity condition on $P$ ensures that $|T(P)|\gg N^3.$ We often write just $T$ for the set $T(P)$.

Let $Q_c(T)$ denote the number of pairs of congruent triangles in $T$ and $Q_s(T)$ the number of pairs of similar triangles. Clearly $Q_s(T)\geq Q_c(T)$, since congruent triangles are similar.

\subsubsection*{Proof of (i)} Given two triangles $\tau_1,\tau_2\in T$, they are congruent if and only if there exists a rigid motion $\phi$ of $\R^2$, such that $\phi(\tau_1)=\tau_2$.  The rigid motion $\phi$ is either homotopic to the identity or not. Let $G_c$ be the group of rigid motions, homotopic to the identity (every such transformation is a composition of a translation and a rotation, but not a mirror symmetry).  Let $\Phi\subset G_c$ be the set of rigid motions, such that for every $\phi\in \Phi$, $\phi(\tau_1)=\tau_2$ for some $\tau_1,\tau_2\in T(P)$. Let $n(\phi)$ denote the number of pairs of congruent triangles, arising from one other via the transformation $\phi$.

It can be assumed without loss of generality that for at least half of the pairs of congruent triangles of $T(P)$, the two triangles are obtained from one another via a rigid motion $\phi$, which is homotopic to the identity. I.e., without loss of generality
\begin{equation}
\sum_{\phi\in \Phi} n(\phi) \;\geq \;\frac{1}{2} Q_c(T).
\label{hi}\end{equation}

Otherwise, the set $\Phi$ would be re-defined as follows. let $P'$ be the mirror image of $P$ against some line, chosen for certainty's sake such that $P$ and $P$ and $P'$ are disjoint. $\Phi$ would now be the set of rigid motions $\phi$, homotopic to the identity, such that for every $\phi\in \Phi$, $\phi(\tau_1)=\tau_2$ for some $\tau_1\in T(P),\tau_2\in T(P').$ It is easy to see that for the rest of the proof this makes no principal difference.

Now, for every pair of points $p,q\in P$, let $L_{pq}$ denote the set of rigid motions in $G_c$, taking $p$ to $q$. It was shown in (\cite{GK}) that the set $L_{pq}$ can be parameterised as a straight line in $\R^3$ (to which one may append the point at infinity corresponding to the translation of $p$ to $q$). As a mater of fact, each $L_{pq}$ is a right coset in $G_c$ of $SO(2)$ by the translation. Indeed, any transformation in the set $L_{pq}$ arises as a composition of the translation of $p$ to $q$ and then rotating around $q$ by an arbitrary angle. It was also shown in \cite{GK} that the set of $N^2$ lines $\{L_{pq}\}_{p,q\in P}$ satisfies the conditions of Theorem \ref{theorem:GKmain}.

A triple intersection of the lines $L_{pq},L_{p'q'}, L_{p''q''}$ would yield $\phi\in T$ and the pair of congruent triangles $pp'p''$ and $qq'q''$ (which may be degenerate, since there are no restrictions on what the three vertices are). Hence, for $k\geq 3$, let $S_k$ be the set of points in $\R^3$, where the number between $k$ and $2k$ distinct lines $L_{pq}$ intersect. By Theorem \ref{theorem:GKmain},
$$
|S_k|\ll \frac{N^3}{k^2}.
$$
Thus, the number of intersecting triples of distinct lines from $\{L_{pq}\}_{p,q\in P}$, and hence the number of pairs of congruent triangles is bounded, via a dyadic summation with $k$ running over the powers $2^j$ of $2$, between $3$ and $N$, as follows:
$$
Q_c(T) \ll \sum_{k = 3, 4,\ldots, 2^j,\ldots, N} k^3 |S_{k}|\ll N^4.
$$
On the other hand, the total number of non-degenerate triangles $|T(P)|\gg N^3$, and therefore, by the Cauchy-Schwarz inequality, the number of distinct congruence classes is
$$
\Omega\left(\frac{N^6}{Q_c(T)} \right)\gg N^2.
$$
This proves the claim (i) of Theorem \ref{mainth}.

\medskip
\subsubsection*{Proof of (ii)}
Let us now view the points of $P$ as complex numbers. Let now $G_s$ denote the group of linear conformal transformations of $\C$, with elements $\psi$:
\begin{equation}\label{cm}\psi(z) = az + b, \mbox{ for some }a,b\in \C,\;a\neq 0.\end{equation}

Let $\Psi\subset G_s$ be the set of linear conformal mappings, such that for every $\psi\in \Psi$, $\psi(\tau_1)=\tau_2$ for some $\tau_1,\tau_2\in T(P)$. That is $\psi$ maps the vertices $p,p',p''$ of $\tau_1$ to the vertices $q,q',q''$ of $\tau_2$, respectively. Let $n(\psi)$ denote the number of similar triangles, arising from one other via the transformation $\psi$.

Just like (\ref{hi}), it can be assumed without loss of generality that
\begin{equation}
\sum_{\psi\in \Psi} n(\psi) \;\geq \;\frac{1}{2} Q_s(T).
\label{hii}\end{equation}
(Otherwise one can once again consider mapping  $P$ to $P'$, where $P'$ is the complex conjugate of $P$ or just use the complex conjugate of $z$ in the definition (\ref{cm}).)

Given $p,q\in P$, let us use the same notation $L_{pq}$ for the set of linear conformal mappings taking $p$ to $q$. The set $L_{pq}$ will now be a line in $\C^2$, as follows:
$$
L_{pq} = \{(a,b)\in \C^2:\,ap + b = q\}
$$
A triple intersection of the lines $L_{pq},L_{p'q'}, L_{p''q''}$ would correspond to now a similarity transformation $\psi$, yielding a pair of similar (homothetic) triangles $pp'p''$ and $qq'q''$.

So, now there is a set of $N^2$ lines in $\C^2$, and no more than $N$ lines are concurrent, for otherwise there will be a linear conformal transformation $\psi$ that takes some $p$ simultaneously to $q$ and $q'\neq q$.

Repeating verbatim what was done in the proof of (i), for $k\geq 3$, let $S_k$ be the set of points in $\C^2$, where the number between $k$ and $2k$ distinct lines $L_{pq}$ intersect. By Theorem \ref{theorem:GKmain},
$$
|S_k|\ll \frac{N^4}{k^3}.
$$
Thus, the number of intersecting triples of distinct lines from $\{L_{pq}\}_{p,q\in P}$, and hence the number of pairs of similar triangles in $T(P)$ is bounded, via a dyadic summation, as follows:
$$
Q_s(T) \ll \sum_{k = 3, 4,\ldots, 2^j,\ldots, N} k^3 |S_{k}|\ll {N^4}\log N.
$$
Then, by the Cauchy-Schwarz inequality, the number of distinct classes of similar triangles is
$$
\Omega\left(\frac{N^6}{Q_s(T)} \right)\gg \frac{N^2}{\log N}.
$$
This completes the proof  of Theorem \ref{mainth}.
\end{proof}


\vskip.5in

\begin{thebibliography}{11}

\vskip.125in



\bibitem{BMP} P. Brass, W.O.J. Moser, J. Pach. {\em Research Problems in Discrete Geometry.} Springer Verlag (2005), 499pp.




\bibitem{E} P. Erd\H os. {\em On sets of distances of $n$ points.} American Mathematical Monthly {\bf 53} (1946), 248-–250.

\bibitem{ES} G. Elekes, M. Sharir. {\it Incidences in three dimensions and distinct distances in the
    plane}. Proceedings 26th ACM Symposium on Computational Geometry (2010), 413--422.

\bibitem{EKS}  G. Elekes, H. Kaplan,  M. Sharir. {\em On lines, joints, and incidences in three
dimensions.} Journal of Combinatorial Theory, Series A, {\bf 118} (2011), 962--977.



\bibitem{GK} L. Guth, N.H. Katz. {\em On the Erd\"os distinct distance problem in the plane.} Preprint {\sf arXiv:math/1011.4105} (2010), 37pp.


\bibitem{IRR}  A. Iosevich, O. Roche-Newton, M. Rudnev. {\em On an application of Guth-Katz theorem. } Math. Res. Lett. {\bf 18} (2011), no. 4, 691–-697.

\bibitem{G}    M. Gromov. {\em Isoperimetry of waists and concentration of maps.} Geom. Funct. Anal. {\bf 13} (2003), 178--215.


\bibitem{KMS} H. Kaplan, J. Matou$\check{\rm s}$ek, M. Sharir. {\em Simple Proofs of Classical Theorems in Discrete Geometry via the Guth--Katz Polynomial Partitioning Technique.} Preprint {\sf arXiv:math/1102.5391} (2011), 17pp.

\bibitem{RR} O Roche-Newton, M. Rudnev. {\em Areas of rectangles and product sets of sum sets.} Preprint {\sf  arXiv:math/1203.6237} (2012), 13pp.



\bibitem{StT} A.H. Stone, J.W. Tukey. {\em Generalized "sandwich'' theorems.} Duke Math. J. {\bf 9} (1942), 356–-359.

\bibitem{ST}   E. Szemer\'edi, W. T. Trotter. {\em Extremal problems in discrete geometry.} Combinatorica {\bf 3} (1983), 381--392.


\bibitem{SoT} J. Solymosi, T. Tao. {\em An Incidence Theorem in Higher Dimensions.}  Discrete and Comp. Geometry  (2012), DOI: 10.1007/s00454-012-9420-x.

\bibitem{SoTa} J. Solymosi, G. Tardos. {\em On the number of k-rich transformations.} Computational geometry (SCG'07), 227–-231, ACM, New York, 2007.

\bibitem{T} C. T\'oth. The Szemer\'edi-Trotter theorem in the complex plane. Preprint {\sf
     arXiv:math/0305283} (2003), 23pp.

\bibitem{Z} J. Zahl. {\em An improved bound on the number of point-surface incidences in three dimensions.}  Preprint {\sf arXiv:math/1104.4987} (2011), 18pp.

\bibitem{Z1} J. Zahl. {\em A Szemer\'edi-Trotter type theorem in $\R^4$.} Preprint {\sf arXiv:math/1203.4600} (2012), 47pp.

\end{thebibliography}
\end{document}